\documentclass[a4paper,10pt,reqno]{amsart}
\usepackage[a4paper,tmargin=40mm,bmargin=35mm,hmargin=30mm]{geometry}
\usepackage[T1]{fontenc}
\usepackage{lmodern}
\usepackage{amsmath}
\usepackage{amssymb}
\usepackage{amsthm}
\usepackage{stmaryrd}
\usepackage{tikz}
\usepackage{booktabs}
\usepackage{multirow}
\usepackage{multicol}
\usepackage{here}
\usepackage{aliascnt}
\usepackage{hyperref}
\usepackage[noabbrev]{cleveref}

\newcommand{\NewTheorem}[2]{
	\newaliascnt{#1}{TheoremEnvironment}
	\newtheorem{#1}[#1]{#1}
	\aliascntresetthe{#1}
	\crefname{#1}{#1}{#2}
	\Crefname{#1}{#1}{#2}
}

\theoremstyle{definition}

\NewTheorem{Definition}{Definitions}
\NewTheorem{Remark}{Remarks}
\NewTheorem{Axiom}{Axioms}
\NewTheorem{Example}{Examples}

\NewTheorem{Observation}{Observations}
\NewTheorem{Convention}{Conventions}
\NewTheorem{Notation}{Notations}
\NewTheorem{Setting}{Settings}
\NewTheorem{Question}{Questions}
\NewTheorem{Answer}{Answers}
\NewTheorem{Conjecture}{Conjectures}
\NewTheorem{Problem}{Problems}
\NewTheorem{Solution}{Solutions}
\NewTheorem{Definition and Remark}{Definition and Remark}
\NewTheorem{Comment}{Comments}
\NewTheorem{Aim}{Aims}
\NewTheorem{Caution}{Cautions}
\NewTheorem{Exercise}{Exercises}

\theoremstyle{plain}
\NewTheorem{Proposition}{Propositions}
\NewTheorem{Lemma}{Lemmas}
\NewTheorem{Theorem}{Theorems}
\NewTheorem{Corollary}{Corollaries}

\crefname{enumi}{}{}
\Crefname{enumi}{}{}
\creflabelformat{enumi}{(#2#1#3)}
\crefname{enumii}{}{}
\Crefname{enumii}{}{}
\creflabelformat{enumii}{(#2#1#3)}
\crefname{enumiii}{}{}
\Crefname{enumiii}{}{}
\creflabelformat{enumiii}{(#2#1#3)}

\makeatletter
\renewcommand{\p@enumii}{}
\renewcommand{\p@enumiii}{}
\makeatother

\numberwithin{equation}{section}
\crefname{equation}{}{}
\Crefname{equation}{}{}
\creflabelformat{equation}{(#2#1#3)}

\newcommand{\SwapSymbols}[1]{
	\expandafter\let\expandafter\temporarysymbol\csname #1\endcsname
	\expandafter\let\csname #1\expandafter\endcsname\csname var#1\endcsname
	\expandafter\let\csname var#1\endcsname\temporarysymbol
}

\SwapSymbols{epsilon}
\SwapSymbols{phi}
\SwapSymbols{Gamma}
\SwapSymbols{Delta}
\SwapSymbols{Theta}
\SwapSymbols{Lambda}
\SwapSymbols{Xi}
\SwapSymbols{Pi}
\SwapSymbols{Sigma}
\SwapSymbols{Upsilon}
\SwapSymbols{Phi}
\SwapSymbols{Psi}
\SwapSymbols{Omega}

\newcommand{\bbZ}{\mathbb{Z}}

\newcommand{\cA}{\mathcal{A}}

\newcommand{\cG}{\mathcal{G}}

\newcommand{\cK}{\mathcal{K}}

\newcommand{\cX}{\mathcal{X}}

\newcommand{\To}{\rightarrow}

\newcommand{\rto}{\rightarrow}

\DeclareMathOperator{\End}{End}

\DeclareMathOperator{\Irr}{IrrSpec}
\DeclareMathOperator{\Mod}{Mod}

\DeclareMathOperator{\Mon}{Mon}

\DeclareMathOperator{\Adim}{Adim}
\DeclareMathOperator{\GKdim}{KG-dim}
\DeclareMathOperator{\kdim}{k-dim}

\DeclareMathOperator{\Kdim}{K-dim}
\DeclareMathOperator{\ckdim}{ck-dim}
\DeclareMathOperator{\Comp}{Comp}

\DeclareMathOperator{\CSupp}{CSupp}

\DeclareMathOperator{\Ann}{Ann}

\DeclareMathOperator{\Ker}{Ker}

\DeclareMathOperator{\Spec}{Spec}

\DeclareMathOperator{\Ass}{Ass}
\DeclareMathOperator{\Supp}{Supp}

\DeclareMathOperator{\ASpec}{ASpec}

\DeclareMathOperator{\AMin}{AMin}
\DeclareMathOperator{\rAnn}{rAnn}
\DeclareMathOperator{\AAss}{AAss}
\DeclareMathOperator{\aass}{aass}
\DeclareMathOperator{\ASupp}{ASupp}

\DeclareMathOperator{\ass}{ass}

\DeclareMathOperator{\MASpec}{m-ASpec}

\DeclareMathOperator{\cm}{cm}

\title{A new dimension for Grothendieck categories via the 
atom spectrum}


\author{Negar Alipour and Reza Sazeedeh}
\address{Department of Mathematics, Urmia University, P.O.Box: 165, Urmia, Iran}
\email{negaralipur8707@yahoo.com}

\address{Department of Mathematics, Urmia University, P.O.Box: 165, Urmia, Iran}
\email{rsazeedeh@ipm.ir and r.sazeedeh@urmia.ac.ir}

\begin{document}

\begin{abstract}
In this paper, we define a new dimension for objects in a Grothendieck category $\cA$. We show that it serves as
 a lower bound for Gabriel-Krull dimension and under certain conditions, the two dimensions coincide. We carry out our investigation for a  fully right bounded ring $A$. We introduce a new spectrum $\Comp A$ via compressible right $A$-modules. In analogy with dimension theory for commutative rings,  we show that the Krull dimension of right $A$-modules can be computed via the length of  chain of  prime  ideals of $A$ and also the length of  chain of elements of $\Comp A$.   
	\end{abstract}

\maketitle

\tableofcontents

\section{Introduction}
The Krull dimension of a topological space is the maximal length of strictly chain of its irreducible closed subsets. Gabriel [G] extended this concept to abelian categories within a purely algebraic framework, where it is known as the Krull-Gabriel dimension. This theory was subsequently generalized to noncommutative rings and representation theory (see [D, MR, GR, GW, MR]). 

From  a different perspective, the spectrum of an abelian category $\cA$ and its associated spectral dimension form a bridge between category theory and geometry.  This concept generalizes the classical notions of prime spectrum of a commutative ring and Krull dimension defined via the length of chain of prime ideals to a much border algebraic context. Gabriel [G] associated to a Grothendieck category a spectrum space  whose points are the isomorphism classes of indecomposable injective objects. This was later generalized by Rosenberg [Ro], who introduced several spectra for an abelian category and raised the following natural question.

{\it What is the relationship between the Krull-Gabriel dimension of an abelian category and its spectral dimension?}

For an abelian category $\cA$, Kanda [K1], introduced the atom spectrum $\ASpec\cA$ equipped with a preorder $\leq $. This construction is inspired by monoform modules and their equivalence relation over noncommutative rings, as explored by Storerr [St]. The main goal of this work is to answer Rosenberg's question where the spectrum is $\ASpec\cA$ and the spectral dimension is defined via $\leq$. Let $A$ be a  fully right bounded ring. We investigate  the spectrum $\Spec A$ and introduce $\Comp A$ consisting of $[M]$, the equivalence class of a compressible right $A$-module $M$. We show that the Krull-Gabriel dimension and the classical Krull dimension of right $A$-modules coincide with the spectral dimensions associated to  these spectra.

Throughout this paper, we assume that $\cA$ is a Grothendieck category.  In Section 2, we study the preorder $\leq$ on $\ASpec\cA$ and we give several preliminaries lemmas related to $\leq$.
In Section 3, we study Krull-Gabriel dimension of objects in $\cA$.
In Section 4, we study spectral dimension of $\ASpec \cA$. For any object $M\in\cA$, we define the atomical dimension  $\Adim M$ via the preorder $\leq $ on $\ASpec \cA$. An object $M$ is said to be locally monoform if any nonzero subquotient of $M$ contains a monoform subobject. The following theorem is the main result of this section.

\begin{Theorem}[\cref{coindim}]
Let $M$ be an object in $\cA$. Then $\Adim M\leq \GKdim M$. Moreover, if $\ASpec\cA$ is Alexandrov and $M$ is locally monoform, then  $\GKdim M=\Adim M$.
\end{Theorem}

In Section 5, we investigate the minimal atoms of an object in $\cA$. For an object $M$ in $\cA$, an atom $\alpha\in\ASupp M$ is said to be minimal if it is minimal in $\ASupp M$ under $\leq$. For any $\alpha\in\ASpec\cA$, we define $\Lambda(\alpha)=\{\beta\in\ASpec\cA|\hspace{0.1cm} \alpha\leq \beta\}$. We show that if $M$ is noetherian and  $\Lambda(\alpha)$ is an open subset of $\ASpec\cA$ for any $\alpha\in\AMin M$, then $\AMin M$ is a finite set; see \cref{finmin}.

 Since the Alexandrov condition on $\ASpec \cA$ plays a crucial role in our investigation, in Section 6, we provide examples of categories with Alexandrov atom spectra. We show that if $A$ is a right fully bounded ring with right Krull dimension, then $\ASpec A$ is Alexandrov; see \cref{two}. We also  show that if $R$ is a commutative ring and $M$ is a noetherian $R$-module with $A=\End_R(M)$, then $\ASpec A$ is Alexandrov; see \cref{comend}.
 
  For a compressible right $A$-module $M$, we denote by $[M]$ the class of compressible modules $N$ such that $\Ann(M)=\Ann(N)$ and we set 
${\rm Comp} A=\{[M]|\hspace{0.05cm}$ $M$ is a compressible right $A$-module$\}$. We prove  the  following theorem.

\begin{Theorem}[\cref{tcorr,bijj,cabij}]
Let $A$ be a right fully bounded ring.  

${\rm (1)}$ If $A$ has a right Krull-dimension, then there is a bijective map $\Ann:\Comp A\To\Spec A$ given by $[M]\mapsto \Ann(M)$. 

Moreover, If $A$ is a right noetherian ring, then

${\rm (2)}$  there is a bijective map $\aass:\Spec A\To\ASpec A$ given by $\frak p\mapsto \aass(\frak p)$; 

${\rm (3)}$  there is a bijective map $\Comp A\To \ASpec A$, given by $[M]\mapsto\overline{\cm(\Ann(M))}$.
\end{Theorem}

 This theorem yields us to define $\Supp M$ and $\CSupp M$, the  support and the compressible support  of a right $A$-module $M$ as a subset of $\Spec A$ and $\Comp A$, respectively. The little  Krull dimension of $M$,  denoted by $\kdim M$, is the supremum of the length of chain of prime ideals in $\Supp M$. Moreover, for a fully right bounded ring $A$, the  compressible Krull dimension of $M$, dented by $\ckdim M$, is the supremum of the length of chain of elements in $\CSupp M$.  We end this paper by the following result. 
\begin{Theorem}[\cref{bkr,last}]
Let $A$ be a fully right bounded ring and $M$ be a right $A$-module with $\Adim M<\omega$. Then $$\GKdim M=\Adim M=\kdim M=\ckdim M\leq \Kdim M.$$ In particular, if $M$ is noetherian, the inequality is equality. 
\end{Theorem}


\section{Atom spectrum of a Grothendieck category}

In this section we recall some basic definitions of category theory, especially    definitions on atom spectrum of a Grothendieck category $\cA$ given by \cite{K2}. 

\begin{Definition}

(1) An abelian category $\cA$ with a generator is said to be a {\it Grothendieck category} if it has arbitrary direct sums and direct limits of short exact sequence are exact, this means that if a direct system of short exact sequences in $\cA$ is given, then the induced sequence of direct limits is a short exact sequence.

(2) A Grothendieck category $\cA$ is said to be {\it locally noetherian} if it has a small generating set of noetherian objects.
\end{Definition}

 \medskip
\begin{Definition}
 A nonzero object $M$ in $\cA$ is {\it monoform} if for any
nonzero subobject $N$ of $M$, there exists no common nonzero
subobject of $M$ and $M/N$ which means that there does not exist a
nonzero subobject of $M$ which is isomorphic to a subobject of
$M/N$. We denote by Mon$\cA$, the set of all monoform
objects in $\cA$.
\end{Definition}

 Two objects $H$ and $H'$ in Mon$\cA$ are said to be {\it
atom-equivalent} if they have a common nonzero subobject. The atom equivalence establishes an
equivalence relation on monoform objects. For every
 object $H$ in Mon$\cA$,  we denoted by $\overline{H}$, the
 {\it equivalence class} of $H$, that is
\begin{center}
  $\overline{H}=\{G\in\Mon\cA|\hspace{0.1cm} H \hspace{0.1cm}
 {\rm and}\hspace{0.1cm} G$ have a common nonzero subobject$\}.$
\end{center}

   \begin{Definition}
 The {\it atom spectrum} $\ASpec\cA$ of $\cA$ is the quotient class
of $\Mon\cA$ consisting of all equivalence classes induced by
this equivalence relation; in other words 
$$\ASpec\cA=\{\overline{H}|\hspace{0.1cm} H\in\Mon\cA\}.$$
Any equivalence class is called an {\it
atom} of $\ASpec\cA$.
\end{Definition}

\medskip
 The notion support and associated prime of  modules over a commutative ring can be generalized to objects in a Grothendieck category $\cA$ via $\ASpec\cA$.  

\begin{Definition}
 Let $M$ be an object in $\cA$.

(1) The {\it atom support} of $M$, denoted by $\ASupp M$, is defined as 
$$\ASupp M=\{\alpha\in\ASpec\cA|\hspace{0.1cm} {\rm there\hspace{0.1cm} exists\hspace{0.1cm}}
 H\in\alpha\hspace{0.1cm} {\rm which \hspace{0.1cm}is\hspace{0.1cm} a
 \hspace{0.1cm}subquotient\hspace{0.1cm} of\hspace{0.1cm}} M\}.$$

(2) The {\it associated atoms} of $M$, denoted by $\AAss M$, is defined as
$$\AAss M=\{\alpha\in\ASupp M|\hspace{0.1cm} {\rm there\hspace{0.1cm} exists\hspace{0.1cm}}
 H\in\alpha\hspace{0.1cm} {\rm which \hspace{0.1cm}is\hspace{0.1cm} a
 \hspace{0.1cm}subobject\hspace{0.1cm} of\hspace{0.1cm}} M\}.$$ 
\end{Definition}

 If $\cG$ is a generating set of objects for $\cA$, any indecomposable injective object in $\cA$ occurs as the injective envelope of some quotient $X/U$ with $X\in\cG$. Then Sp$\cA$, the the collection of  isomorphism classes of indecomposable injective objects in $\cA$ form a set. It is straightforward to show that $E(\alpha)$ is an indecomposable injective object of $\cA$ so that the injective map $\alpha\mapsto E(\alpha)$ establishes $\ASpec\cA$ as a subset of Sp$\cA$. Open subsets of $\ASpec\cA$ can be defined in correspondence with specialization closed subsets of $\Spec A$, where $A$ is a commutative ring.

\begin{Definition}
A subset $\Phi$ of $\ASpec\cA$ is said to be {\it open} if for any
$\alpha\in\Phi$, there exists a monoform  object $H$ in $\cA$ such that $\alpha=\overline{H}$ and $\ASupp H\subset\Phi$. For every nonzero object $M$ in $\cA$, it is clear that $\ASupp M$ is an open subset of $\ASpec\cA$. 
 \end{Definition}

 We recall from [K1] that $\ASpec\cA$ can be regarded as a preordered set together with a specialization order $\leq$ as follows:
for any atoms $\alpha$ and $\beta$ in $\ASpec\cA$, we have 
$\alpha\leq \beta$ if and only if for any open subset $\Phi$ of
$\ASpec\cA$ satisfying $\alpha\in \Phi$, we have $\beta\in\Phi$.

\medskip

\begin{Definition and Remark}\label{def}
 An atom $\alpha$ in $\ASpec\cA$ is said to be {\it maximal} if there
exists a simple object $H$ in $\cA$ such that $\alpha=\overline{H}$. The class of all maximal atoms in $\ASpec\cA$ is denoted by $\MASpec \cA$. If $\alpha$ is a maximal atom, then
$\alpha$ is maximal in $\ASpec\cA$ under the order $\leq$; see [S1, Remark 4.7].
\end{Definition and Remark}
  \medskip
 \begin{Definition}
	A topological space $X$ is said to be {\it Alexandrov} if the intersection of any family of open subsets of $X$ is open. 
	\end{Definition} 
  
	\medskip
\begin{Definition}
A full subcategory $\cX$ of an abelian category $\cA$ is said to be {\it Serre} if for any exact sequence $0\To M\To N\To K\To 0$ in $\cA$, the object $N$ belongs to $\cX$ if and only if $M$ and $K$ belong to $\cX$. A Serre subcategory $\cX$ of $\cA$ is said to be {\it localizing} if the
canonical functor $F:\cA\To\cA/\cX$ has a right adjoint functor
$G:\cA/\cX\To \cA$. 
\end{Definition}

\medskip

For every $\alpha\in\ASpec\cA$, the topological closure of $\alpha$, denoted by $\overline{\{\alpha\}}$ consists of all
 $\beta\in\ASpec\cA$ such that $\beta\leq\alpha$. The localizing subcategory induced by $\alpha$ is 
 $\cX({\alpha})=\ASupp^{-1}(\ASpec\cA\setminus\overline{\{\alpha\}})$.

\medskip
Let $\cX$ be a localizing subcategory of $\cA$ and $\alpha\in\ASpec\cA\setminus\ASupp \cX$. Then for any monoform $H$ in $\cA$ with $\overline{H}=\alpha$, we have $H\notin\cX$ and so it follows from [K1, Lemma 5.14] that $F(H)$ is a monoform object of $\cA/\cX$ where $F:\cA\rto \cA/\cX$ is the canonical functor. In this case, we denote $\overline{F(H)}$ by $F(\alpha)$. 

Suppose that $\alpha\in\ASpec\cA/\cX$ and $H_1, H_2$ are monoform objects in $\cA/\cX$ with $\alpha=\overline{H_1}=\overline{H_2}$. Then  $G(H_1)$ and $G(H_2)$ are monoform objects in $\cA$ by [K1, Lemma 5.14] and since $G$ is faithful, they have a common nonzero subobject. We denote $\overline{G(H_1)}=\overline{G(H_2)}$ by $G(\alpha)$. Hence we have functions $F:\ASpec\cA\setminus \ASupp\cX\To \ASpec\cA/\cX$ given by $\alpha\mapsto F(\alpha)$ and $G:\ASpec\cA/\cX\To \ASpec\cA\setminus\ASupp\cX$ given by $\alpha\mapsto G(\alpha)$. In the following lemma we show that there is a bijection between $\ASpec\cA/\cX$ and $\ASpec\cA\setminus\ASupp\cX$.
\medskip

\begin{Lemma}\label{bij}
The map $F:\ASpec\cA\setminus \ASupp\cX\To \ASpec\cA/\cX$ is a homeomorphism with the inverse map $G$.  
\end{Lemma}
\begin{proof}
See [K1, Theorem 5.17].
\end{proof}

\begin{Lemma}\label{leqf}
If  $\alpha_1,\alpha_2\in\ASpec\cA\setminus\ASupp\cX$  such that $\alpha_1\leq\alpha_2$, then $F(\alpha_1)\leq F(\alpha_2)$
\end{Lemma}
\begin{proof}
It is clear that $\{\ASupp F(M)|\hspace{0.1cm} M\in\cA\}$ forms a basis of open subsets of $\ASpec\cA/\cX$ and $F(\ASupp M\setminus\ASupp\cX)=\ASupp F(M)$ for any object $M$ in $\cA$. Now the result follows by the definition of $\leq$
\end{proof}

A similar argument yields the following lemma.
\begin{Lemma}\label{leq}
Let $\alpha_1,\alpha_2\in\ASpec\cA/\cX$ such that $\alpha_1\leq \alpha_2$. Then $G(\alpha_1)\leq G(\alpha_2)$.
\end{Lemma}
  
  \medskip
For any $\alpha\in\ASpec\cA$, we define $\Lambda(\alpha)=\{\beta\in\ASpec\cA|\hspace{0.1cm} \alpha\leq \beta\}$ and we have the following lemma.

 \begin{Lemma}\label{SaS}
 Let $\alpha\in\ASpec\cA$. Then $\Lambda(\alpha)=\bigcap_{\overline{H}=\alpha}\ASupp H$.
 \end{Lemma}
 \begin{proof}
 See \cite[Proposition 2.3]{SaS}.
 \end{proof}

\begin{Lemma}\label{lamop}
For every atom $\alpha\in\ASpec\cA\setminus\ASupp\cX$, we have $F(\Lambda(\alpha))=\Lambda(F(\alpha))$. Moreover, if $\Lambda(\alpha)$ is an open subset of $\ASpec\cA$, then $\Lambda(F(\alpha))$ is open in $\ASpec\cA/\cX$.
\end{Lemma}
\begin{proof}
The first assertion is straightforward by using \cref{bij}, \cref{leq} and \cref{leqf}. Given $F(\beta)\in\Lambda(F(\alpha))$, according to \cref{bij} and \cref{leq}, we have $\alpha<\beta$ so that $\beta\in\Lambda(\alpha)$. Then by the assumption, there exists a monoform object $H$ in $\cA$ such that $\beta=\overline{H}$ and $\ASupp H\subseteq\Lambda(\alpha)$. It follows from \cref{bij} and the first assertion that $\ASupp F(H)=F(\ASupp H\setminus\ASupp\cX)\subseteq\Lambda(F(\alpha))$. Therefore the result follows as $F(\beta)=\overline{F(H)}$. 
\end{proof}

 \medskip
The following result holds for a general topology as well.
\medskip
\begin{Lemma}\label{alex}
$\ASpec\cA$ is an Alexandrov topological space if and only if $\Lambda(\alpha)$ is an open subset of $\ASpec \cA$ for all $\alpha\in\ASpec\cA$.
\end{Lemma}
\begin{proof}
The part "only if"' follows from \cref{SaS}. Assume that $\{{U_i}\}_{i\in\Gamma}$ is a family of open subsets of $\ASpec\cA$ and $\alpha\in\bigcap_{\Gamma}\limits U_i$ is an arbitrary atom. As $\alpha\in U_i$ for each $i$, there is a monoform $H_i$ such that $\alpha=\overline{H_i}$ and $\ASupp H_i\subseteq U_i$ for each $i$. Since $\Lambda(\alpha)$ is open, using again \cref{SaS}, we have $\alpha\in\Lambda(\alpha)=\bigcap_{\overline{H}=\alpha}\limits\ASupp H\subseteq\bigcap_{\Gamma}\limits\ASupp H_i\subseteq\bigcap_{\Gamma}\limits U_i$.
\end{proof}

\begin{Lemma}\label{lx}
If $\ASpec\cA$ is Alexandrov, then so is $\ASpec\cA/\cX$.
\end{Lemma}
\begin{proof}
It is straightforward by using \cref{alex} and \cref{lamop}. 
\end{proof}

\begin{Definition}
An object $M$ in $\cA$ is {\it locally monoform} if any subquotient of $M$ contains a monoform subobject. An atom $\alpha\in\ASpec \cA$ is locally monoform if there exists a monoform locally monoform object $H$ such that $\alpha=\overline{H}$. The category $\cA$ is said to be locally monoform if any nonzero object in $\cA$ is locally monoform.  
\end{Definition}

\medskip
\begin{Lemma}\label{alexmax}
Let $\alpha$ be a locally monoform atom in $\ASpec\cA$ such that $\Lambda(\alpha)$ is an open subset of $\ASpec\cA$. Then $\alpha$ is maximal if and only if it is maximal under $\leq$.
\end{Lemma}
\begin{proof}
If $\alpha$ is a maximal atom, in view of \cref{def}, it is maximal under $\leq$.  Now, assume that $\alpha\in\ASpec\cA$ is maximal under $\leq$. According to \cref{SaS} and the assumption, $\Lambda(\alpha)=\{\alpha\}$ is an open subset of $\ASpec \cA$ and so there exists a monoform object $H$ in $\cA$ such that $\ASupp H=\{\alpha\}$ and $\alpha=\overline{H}$. Since $\alpha$ is locally monoform and any nonzero subobject of a locally monoform object is locally monoform, we may assume that $H$ is locally monoform. If $H$ is not simple, it contains a proper  subobject $N$ and so $H/N$ contains a monoform subobject $M$. But since $\ASupp H/N\subseteq \{\alpha\}$, we deduce that $\alpha=\overline{M}$ which is a contradiction.     
\end{proof}

\section{Krull-Gabriel dimension of objects}

Let us begin this section with a definition.

\begin{Definition}\label{Gab}
For a Grothendieck category $\cA$, we define the {\it Krull-Gabriel filtration} of $\cA$ as follows. For any ordinal (i.e ordinal number) $\sigma$
we denote by $\cA_{\sigma}$, the localizing subcategory of $\cA$ which is defined in the following manner:

$\cA_{-1}$ is the zero subcategory.

$\cA_0$ is the smallest localizing subcategory containing all simple objects.

Let us assume that $\sigma=\rho+1$ and denote by $F_{\rho}:\cA\To\cA/\cA_{\rho}$ the canonical functor and by 
$G_{\rho}:\cA/\cA_{\rho}\To\cA$ the right adjoint functor of $F_{\rho}$. Then an object $X$ in $\cA$ will belong to 
$\cA_{\sigma}$ if and only if $F_{\rho}(X)\in{\rm Ob}(\cA/\cA_{\rho})_0$. The left exact radical functor (torsion functor) corresponding to $\cA_{\rho}$ is denoted by $t_{\rho}$.
If $\sigma$ is a limit ordinal, then $\cA_{\sigma}$ is the localizing subcategory generated by all localizing subcategories $\cA_{\rho}$ with $\rho< \sigma$. It is clear that if $\sigma\leq \sigma'$, then $\cA_{\sigma}\subseteq\cA_{\sigma'}$. Moreover, since the class of all localizing subcategories of $\cA$ is a set, there exists an ordinal $\tau$ such that $\cA_{\sigma}=\cA_{\tau}$ for all $\sigma \leq \tau$. Let us put $\cA_{\tau}=\cup_{\sigma}\cA_{\sigma}$. Then $\cA$ is said to be {\it semi-noetherian}  if $\cA=\cA_{\tau}$. We also say that the localizing subcategories $\{\cA_{\sigma}\}_{\sigma}$ define the Krull-Gabriel filtration of $\cA$. We say that an object $M$ in $\cA$ has the {\it Krull-Gabriel dimension} defined or $M$ is {\it semi-noetherian} if $M\in{\rm Ob}(\cA_{\tau})$. The smallest ordinal  $\sigma$ so that $M\in{\rm Ob}(\cA_{\sigma})$ is denoted by $\GKdim M$. Because the class of ordinals is well-ordered, throughout this paper, $\omega$ is denoted the smallest limit ordinal. We observe that $\GKdim 0=-1$ and $\GKdim M\leq 0$ if and only if $\ASupp M\subseteq \textnormal{m-}\ASpec \cA$. 
\end{Definition}

\begin{Definition}
For any atom $\alpha\in\ASpec \cA$, the {\it Gabriel-Krull dimension of $\alpha$}, is the least ordinal $\sigma$ such that $\alpha\in\ASupp\cA_\sigma$. If such an ordinal exists, we denote it by $\GKdim \alpha$ and by definition, it is a non-limit ordinal.
\end{Definition}

\medskip
	\begin{Lemma}\label{dist}
	Let $\alpha,\beta$ be two atoms in $\ASupp\cA$ such that $\GKdim\alpha=\GKdim\beta$. Then $\alpha\nless\beta$.
	\end{Lemma}
	\begin{proof}
If $\alpha<\beta$, then $\sigma$ is a non-limit ordinal by the definition. Using \cref{leqf}, we have $F_{\sigma-1}(\alpha)< F_{\sigma-1}(\beta)$. Since $\alpha,\beta\in\ASupp\cA_{\sigma}$, we have $F_{\sigma-1}(\alpha),F_{\sigma-1}(\beta)\in\ASupp F_{\sigma-1}(\cA_{\sigma})=\ASupp(\cA/\cA_{\sigma})_0$ so that $F_{\sigma-1}(\alpha)$ and $F_{\sigma-1}(\beta)$ are maximal. Thus $F_{\sigma-1}(\alpha)=F_{\sigma-1}(\beta)$; and consequently $\alpha=G_{\sigma-1} F_{\sigma-1}(\alpha)=G_{\sigma-1}F_{\sigma-1}(\beta)=\beta$ by \cref{bij} which is a contradiction.
	\end{proof}

	\medskip
\begin{Corollary}\label{corleq}
If $\alpha,\beta$ are two atoms in $\ASpec\cA$ such that $\alpha<\beta$, Then $\GKdim\beta<\GKdim\alpha$.  
\end{Corollary}
\begin{proof}
Assume that $\GKdim\alpha=\sigma$ for some ordinal $\sigma$. Since $\ASupp\cA_{\sigma}$ is an open subset of $\ASpec\cA$, by definition $\beta\in\ASupp\cA_{\sigma}$. Therefore $\GKdim\beta\leq \sigma$. Now \cref{dist} implies that $\GKdim\beta<\sigma$.
\end{proof}

The following lemma plays a crucial role in our subsequent investigation.

\begin{Lemma}\label{equicat}
For any ordinal $\sigma$, there is $F_0(\cA_\sigma)\cong\cA_\sigma/\cA_0.$ Moreover we have

$$F_0(\cA_\sigma)=\begin{cases}
(\cA/\cA_0)_{\sigma-1} \hspace{0.4cm}{\rm if  \hspace{0.1cm}\sigma<\omega}\\
  (\cA/\cA_0)_{\sigma} 
 \hspace{0.5cm}{\rm if} \hspace{0.1cm}\sigma\geq \omega.
 \end{cases}$$
 \end{Lemma}
\begin{proof}
The equivalence follows from the construction of $\cA/\cA_0$ (see [G, p. 365]) and so it suffices to prove the equalities. We proceed by induction on $\sigma$. If $\sigma<\omega$, then the cases $\sigma=0,1$ are clear  by the definition. Assume that $\sigma>1$ and so by the induction hypothesis, there are the equivalence and equality of categories $$\cA/\cA_{\sigma-1}\cong \cA/\cA_0/F_0(\cA_{\sigma-1})=\cA/\cA_0/(\cA/\cA_0)_{\sigma-2}.$$ If $F'_{\sigma-2}:\cA/\cA_0\to \cA/\cA_0/(\cA/\cA_0)_{\sigma-2}$ is the canonical functor, it suffices to show that $F'_{\sigma-2}(F_0(\cA_\sigma))$ is the smallest subcategory of $\cA/\cA_0/(\cA/\cA_0)_{\sigma-2}$ generated by simple objects. Suppose that $\theta:\cA/\cA_0/(\cA/\cA_0)_{\sigma-2}\To \cA/\cA_{\sigma-1}$ is the equivalence functor.
We now show that $\Ker(\theta \circ F'_{\sigma-2}\circ F_0)=\cA_{\sigma-1}.$ The inclusion $\cA_{\sigma-1}\subset\Ker(\theta \circ F'_{\sigma-2}\circ F_0)$ is clear. If $M\in\Ker(\theta \circ F'_{\sigma-2}\circ F_0)$, then $F'_{\sigma-2}\circ F_0(M)=0$. If $M\notin\cA_{\sigma-1}$, the $M/t_{\sigma-1}(M)$ is nonzero and $F'_{\sigma-2}\circ F_0(M/t_{\sigma-1}(M))=0$. Then by the induction hypothesis, there exists an object $N\in\cA_{\sigma-1}$ such that $F_0(M/t_{\sigma-1}(M))=F_0(N)$. We also may assume that $t_0(N)=0$. Now since $N$ is an essential subobject of $G_0(F_0(N))$ and $M/t_{\sigma-1}(M)$ is a subobject $G_0(F_0(N))$, the objects $M/t_{\sigma-1}(M)$ and $N$ have a common nonzero subobject which is a contradiction. Hence by \cite[Chap 4, p. 180, Theorem 4.9]{Po}, we may assume that $F_{\sigma-1}=\theta \circ F'_{\sigma-2}\circ F_0$; and hence there are the following equalities and equivalences of categories 
$$(\cA/\cA_{\sigma-1})_0=F_{\sigma-1}(\cA_\sigma)=\theta(F'_{\sigma-2}(F_0(\cA_\sigma)))\cong F'_{\sigma-2}(F_0(\cA_\sigma))$$ 
which proves the assertion. We now prove the case $\sigma\geq \omega$. If $\sigma $ is limit ordinal, then $\cA_{\sigma}=\langle\cup_{\rho<\sigma}\cA_{\rho}\rangle_{\rm loc}$ is the localizing subcategory of $\cA$ generated by $\cup_{\rho<\sigma}\cA_{\rho}$ and since $F_0$ is exact and preserves arbitrary direct sums, the induction hypothesis yields $F_0(\cA_{\sigma})= \langle\cup_{\rho<\sigma}(\cA/\cA_0)_{\rho}\rangle_{\rm loc}=(\cA/\cA_0)_{\sigma}$. If $\sigma$ is a non-limit ordinal, then $F_{\sigma-1}$ can be factored as $\cA\stackrel{F_0}\To \cA/\cA_0 \stackrel{F'_{\sigma-1}}\To \cA/\cA_0/(\cA/\cA_0)_{\sigma-1}\cong \cA/\cA_{\sigma-1}.$  Thus $F_{\sigma-1}(\cA_{\sigma})=F'_{\sigma-1}(F_0(\cA_{\sigma}))$ so that $F_0(\cA_{\sigma})=(\cA/\cA_0)_{\sigma}$.
\end{proof}

\medskip
\begin{Lemma}\label{lemcru}
Let $\cX$ be a localizing subcategory of $\cA$ and let $M$ be an  object in $\cA$ such that $\ASupp M\subset\ASupp \cX$. Consider the following conditions.

${\rm (i)}$  $M$ has Krull-Gabriel dimension,

${\rm (ii)}$ $M$ is locally monoform'

 ${\rm (ii)}$ $M\in\cX$.\\Then $(i)\Longrightarrow (ii)\Longrightarrow (iii)$ holds.

  \end{Lemma}
\begin{proof}
See [S2, Lemma 2.10 and Corollary 2.18].
\end{proof}

\medskip

\begin{Lemma}\label{lmono}
\ Let $\alpha\in\ASpec\cA$.

{\rm (i)} If  $\GKdim \alpha$ exists, then $\alpha$ is locally monoform.

{\rm (ii)} If $\alpha$ is locally monoform and $\beta$ is another 
atom in $\ASpec\cA$ such that $\alpha<\beta$, then $\beta$ is locally monoform.

{\rm (iii)} If $M$ is a locally monoform object in $\cA$, then any atom $\alpha\in\ASupp M$ is locally monoform.

{\rm (iv)} If $\cX$ is a localizing subcategory of $\cA$ with the canonical functor $F:\cA\To \cA/\cX$ and $M$ is a locally monoform object in $\cA$, then $F(M)$ is a locally monoform object in $\cA/\cX$. In particular $F(\alpha)$ is locally monoform if $\alpha\in\ASpec\cA\setminus\ASupp\cX$.
\end{Lemma}
\begin{proof}
(i) There exists a monoform object $H$ with Krull-Gabriel dimension such that $\alpha=\overline{H}$. Now, according to \cref{lemcru}, $H$ is locally monoform. (ii) There exits a monoform and locally monoform object $H$ in $\cA$ such that $\alpha\in\ASupp H$. Then $\beta\in\ASupp H$ and so there exists a monoform object $G$ in $\cA$ such that $\beta=\overline{G}$ and $G$ is a subquotient of $H$. This implies that $G$ is locally monoform. (iii) is clear. (iv) If $T$ is a subquotient of $F(M)$, then there exists a subobject $N$ of $M$ and a subquotient $L$ of $M/N$ such that $F(L)=T$. We may assume that $t_\cX(L)=0$ and since $M$ is locally monoform, $L$ contains a monoform subobject $K$. This implies that $F(K)$ is a monoform subobject of $T$. The second assertion follows by the first part.
\end{proof}

\begin{Corollary}\label{coroo}
Let $\sigma$ be an ordinal and let $M$ be an object in $\cA$ such that $\GKdim M=\sigma$. Then $$\GKdim M=\begin{cases}
\GKdim F_0(M)+1\hspace{0.4cm}{\rm if  \hspace{0.1cm}\sigma<\omega}\\
  \GKdim F_0(M) 
 \hspace{0.5cm}{\rm if} \hspace{0.1cm}\sigma\geq \omega 
 \end{cases}$$
\end{Corollary}
\begin{proof}
Assume that $\GKdim M=\sigma$ for some  ordinal $\sigma$. If $\sigma<\omega$, then we have $M\in\cA_\sigma$ and so \cref{equicat} implies that $F_0(M)\in(\cA/\cA_0)_{\sigma-1}$ so that $\GKdim F_0(M)\leq \sigma-1$. If $\GKdim F_0(M)=\rho<\sigma-1$, Then $F_0(M)\in(\cA/\cA_0)_\rho=F_0(\cA_{\rho+1})$ by \cref{equicat}. Thus, there exists an object $N\in\cA_{\rho+1}$ such that $F_0(M)= F_0(N)$. For every $\alpha\in\ASupp M$, if $\alpha$ is a maximal atom, then $\alpha\in\ASupp\cA_0\subseteq \ASupp\cA_{\rho+1}$. If $\alpha$ is not maximal, then $F_0(\alpha)\in\ASupp F_0(M)=\ASupp F_0(N)=F_0(\ASupp N\setminus\ASupp\cA_0)$ so that $\alpha\in\ASupp N$ by \cref{bij}. Hence $\ASupp M\subseteq\ASupp\cA_{\rho+1}$ and so  $M\in\cA_{\rho+1}$ by \cref{lemcru}. Therefore $\GKdim M\leq \rho+1<\sigma$ which is a contradiction. If $\sigma\geq \omega$, it follows from \cref{equicat} that $F_0(M)\in(\cA/\cA_0)_{\sigma}$ and so $\GKdim F_0(M)\leq \sigma$. If $\GKdim F_0(M)=\rho<\sigma$, it follows from the first case that $\rho\geq \omega$. Thus according to \cref{equicat}, we have  $F_0(M)\in(\cA/\cA_0)_{\rho}=F_0(\cA_{\rho})$ and so $\GKdim M\leq \rho<\sigma$ which is a contradiction. 
\end{proof}

\medskip

\begin{Definition}
	Given an ordinal $\sigma\geq 0$, we recall from [MR or GW] that an object $M$ in $\cA$ is  $\sigma$-{\it critical} provided $\GKdim M=\sigma$ while $\GKdim M/N<\sigma$ for all nonzero subobjects $N$ of $M$. It is clear that any nonzero subobject of a $\sigma$-critical object is $\sigma$-critical. An object $M$ is said to be {\it critical} if it is $\sigma$-critical for some ordinal $\sigma$. It is clear to see that any critical object is monoform.  
	\end{Definition}

\medskip
Let $\sigma\leq \GKdim\cA$ be an ordinal number. We denote by $\ASpec_{\sigma}\cA$ the subset of $\ASpec\cA$ consisting of all $\alpha\in\ASpec\cA$ such that $\GKdim\alpha =\sigma+1$. It is clear that for $\sigma\neq\lambda$, we have $\ASpec_{\sigma}\cA\cap \ASpec_{\lambda}\cA=\emptyset$ and if $\cA$ is a semi-noetherian category, then  $\ASpec\cA=\bigcup_{\sigma}\ASpec_{\sigma}\cA$.  For any $\alpha\in\ASpec_{\sigma}\cA$, there exists a monoform object $H$ in $\cA_{\sigma+1}$ such that $\alpha=\overline{H}$. It now follows from \cite[Lemma 2.17]{S2} that $H$ contains a $\sigma+1$-critical subobject $C(\alpha)$. It is clear that $t_{\sigma}(C(\alpha))=0$.

\medskip
\begin{Proposition}\label{p1}
Let $\sigma
$ be an ordinal. Then the map  $\alpha\mapsto F_{\sigma}(C(\alpha))$ provides a bijection between $\ASpec_{\sigma}\cA$ and the set of isomorphism class of simple objects in $\cA/\cA_{\sigma}$. The inverse map is given by $S\mapsto \overline{G_{\sigma}(S)}$. 
\end{Proposition}
\begin{proof}
 Clearly, $F_{\sigma}(C(\alpha))$ is an simple object in $\cA/\cA_{\sigma}$. Hence,  
$\overline{G_{\sigma}(F_{\sigma}(C(\alpha)))}=\alpha$ as $C(\alpha)$ is a subobject of $G_{\sigma}(F_{\sigma}(C(\alpha))).$ Furthermore, if $S$ is a simple object in $\cA/\cA_\sigma$, by [S2, Lemma 2.15], there exists a $\sigma+1$-critical object $M\in\cA$ such that $F_\sigma(M)=S$.  Then $$F_\sigma(\overline{G_\sigma(S)})=
F_\sigma(\overline{G_\sigma(F_\sigma(M))})=
F_\sigma(\overline{M})=F_\sigma(M)=S.$$  
\end{proof}


\section{A new dimension for objects via the atom spectrum}

We begin this section by introducing a new dimension for objects in $\cA$ via the preorder $\leq$ on $\ASpec \cA$.

\begin{Definition}

For every $\alpha\in\ASpec\cA$, we define $\Adim \alpha$, the {\it atomical dimension of} $\alpha$ by transfinite induction.  We say that $\Adim \alpha=0$ if $\alpha$ is maximal under $\leq $. For an ordinal $\sigma>0$, we say that $\Adim\alpha\leq \sigma$ if for every $\beta\in\ASpec\cA$ with $\alpha<\beta$, we have dim $\beta<\sigma$. The least such an ordinal $\sigma$ is said to be the atomical dimension of $\alpha$ and we say that $\Adim \alpha=\sigma$. We set $\Adim 0=-1$. If $\Adim\alpha=n$ is finite, then there exists a chain of atoms $\alpha<\alpha_1<\dots<\alpha_n$ in $\ASpec\cA$ and this chain has the largest length among those starting with $\alpha$.  
For every object $M$ in $\cA$, the {\it atomical dimension} of $M$, denoted by $\Adim M$, is the supremum of $\Adim\alpha$ such that $\alpha\in\ASupp M$.
\end{Definition}

\medskip
\begin{Lemma}\label{mondim}
Let $\alpha$ be an atom in $\ASpec\cA$ such that $\Lambda(\alpha)$ is an open subset of $\ASpec\cA$. Then there exists a monoform object $M$ in $\cA$ such that $\alpha=\overline{M}$ and $\Adim M=\Adim\alpha$. 
\end{Lemma}
\begin{proof}
Since $\Lambda(\alpha)$ is an open subset of $\ASpec\cA$, there exists a monoform object $M$ in $\cA$ such that $\alpha=\overline{M}$ and $\Supp M=\Lambda(\alpha)$. Therefore $\Adim M=\Adim\alpha$. 
 \end{proof}

\medskip

\begin{Lemma}\label{dimm}
Let $\alpha\in\ASpec\cA$. Then we have the following inequalities $$\Adim\alpha\leq \begin{cases}
\Adim F_0(\alpha)+1\hspace{0.5cm}{\rm if} \hspace{0.1cm}\Adim \alpha< \omega\\
  \Adim F_0(\alpha) 
 \hspace{0.5cm}{\rm if} \hspace{0.1cm}\Adim\alpha\geq \omega.
 \end{cases}$$
 Moreover, if  $\ASpec\cA$ is Alexandrov and $\alpha$ is locally monoform, then the inequalities are equalities. 
\end{Lemma}
\begin{proof}
We proceed by transfinite induction on $\Adim \alpha=\sigma$. We first assume that $\sigma<\omega$. The case $\sigma=0$ is clear. If $\sigma>0$, there exists an atom $\beta\in\ASpec\cA$ such that
$\alpha<\beta$ and $\Adim \beta=\sigma-1$. The induction hypothesis implies that $\Adim F_0(\beta)\geq\sigma-2$ so that $\Adim F_0(\alpha)\geq \sigma-1$. To prove the second claim in this case, if $\sigma=0$, by \cref{alexmax}, the atom $\alpha$ is maximal and so there exists a simple object $S$ in $\cA$ such that $\alpha=\overline{S}$. Then $F_0(S)=0$ and  so $\Adim F_0(\alpha)=-1$ by the definition. If $\sigma>0$ and dim$F_0(\alpha)>\sigma-1$, there exists $\beta\in\ASpec\cA$ such that $F_0(\alpha)<F_0(\beta)$ and $\Adim F_0(\beta)\geq\sigma-1$. But \cref{leq} and \cref{bij} imply that $\alpha<\beta$ and the induction hypothesis and \cref{lmono} imply that $\Adim\beta\geq\sigma$ which is a contradiction.  We now assume that $\sigma\geq \omega$. If $\sigma=\omega$, then for any non-negative integer $n$ there exists $\beta\in\ASpec\cA$ such that $\alpha<\beta$ and $\Adim \beta\geq n+1$ and so the first case implies that $\Adim F_0(\beta)\geq n$ so that $\Adim F_0(\alpha)\geq \omega$. Now, assume that $\sigma>\omega$. If $\sigma$ is a non-limit ordinal, then there exists $\beta\in\ASpec\cA$ such that $\alpha<\beta$ and  $\Adim\beta=\sigma-1$.  Thus the induction hypothesis implies that $\Adim F_0(\beta)\geq \sigma -1$ and consequently $\Adim F_0(\alpha)\geq \sigma$. If $\sigma$ is a limit ordinal, it suffices to prove the claim for the case $\sigma=\omega$. For every integer $\rho<\sigma$ there exists $\beta\in\ASpec\cA$ such that $\alpha<\beta $ and $\Adim\beta\geq \rho+1$. The induction hypothesis implies that $\Adim F_0(\beta)\geq \rho$ so that $\Adim F_0(\alpha)\geq \sigma$.
 To prove the second claim in this case, assume that $\ASpec\cA$ is Alexandrov and $\sigma=\omega$. Then for every $\beta\in\ASpec\cA\setminus\ASpec\cA_0$ with $\alpha<\beta$, we have $\Adim \beta<\omega$. Then using the first case, $\Adim F_0(\beta)<\omega$ and hence $\Adim F_0(\alpha)=\omega$. If $\sigma>\omega$ and $\Adim F_0(\alpha)>\sigma$, then there exists  $\beta\in\ASpec\cA\setminus\ASpec\cA_0$ with $\alpha<\beta$ and $\Adim F_0(\beta)\geq \sigma$. But the induction hypothesis  and \cref{lmono} imply that $\Adim \beta=\Adim F_0(\beta)\geq \sigma$ which is a contradiction. 
\end{proof}

\medskip

\begin{Corollary}\label{coroe}
Let $M$ be an object in $\cA$. Then we have the following inequalities 
$$\Adim M\geq \begin{cases}
\Adim F_0(M)+1\hspace{0.5cm}{\rm if} \hspace{0.1cm}\Adim M< \omega\\
  \Adim M 
 \hspace{0.5cm}{\rm if} \hspace{0.1cm}\Adim F_0(M)\geq \omega.
 \end{cases}$$
 Moreover, if $\ASpec\cA$ is Alexandrov and $M$ is locally monoform, then the inequalities are equalities.
\end{Corollary}
\begin{proof}
Straightforward using \cref{dimm}. 
\end{proof}

 The following theorem is the main result of this section showing that the atomical dimension of an object serves as a lower bound for its  Krull-Gabriel dimension.
 \medskip
\begin{Theorem}\label{coindim}
Let $M$ be an object in $\cA$. Then $\Adim M\leq \GKdim M$. Moreover, if $\ASpec\cA$ is Alexandrov and $M$ is locally monoform, then  $\GKdim M=\Adim M$.
\end{Theorem}
\begin{proof}
Assume that $\GKdim M=\sigma$ for some ordinal $\sigma$. We proceed by transfinite induction on $\sigma$. If $\sigma=0$, then $M\in\cA_0$ and so using \cref{def}, every atom in $\ASupp M$ is maximal. Therefore, every atom in $\ASupp M$ is maximal under $\leq$ so that $\Adim M=0$. Suppose inductively that $\sigma>0$ and $\alpha$ is an arbitrary atom in $\ASupp M$. We prove that $\Adim\alpha\leq \sigma$; and consequently $\Adim M\leq \sigma$. For every $\beta\in\ASpec\cA$ with $\alpha<\beta$ and $\GKdim\beta=\rho$, according to \cref{corleq}, we have $\rho<\GKdim \alpha\leq\sigma$. Since $\beta\in\cA_\rho$, there exists $X\in \cA_\rho$ and a monoform object $G$ in $\cA$ such that $\beta=\overline{G}$ and $G$ is a subquotient of $X$. This implies that $\GKdim G\leq \rho$. Now, the induction hypothesis implies that $\Adim\beta\leq \Adim G\leq \GKdim G=\rho<\sigma$; and consequently $\Adim \alpha\leq \sigma$ by the definition.

 To prove the equality, we proceed  by induction on $\Adim M=\lambda$ that $\GKdim M\leq\Adim M$. For the case $\lambda=0$, since $M$ is locally monoform, by \cref{lmono}, any atom in $\ASupp M$ is locally monoform; and hence \cref{alexmax} implies that $\ASupp M\subseteq\ASupp\cA_0$. Therefore, $M\in\cA_0$ by \cref{lemcru} so that $\GKdim M=0$.  Assume that $\lambda>0$. In this case for every $\alpha\in\ASupp M$, by \cref{mondim}, there exists a monoform object $H$ in $\cA$ such that $\alpha=\overline{H}$ and $\Adim\alpha=\Adim H$. If $\Adim\alpha<\lambda$, the induction hypothesis implies that $H\in\cA_{\lambda}$ so that $\alpha\in\ASupp\cA_{\lambda}$. If $\Adim \alpha=\lambda$, we investigate two cases: 
  {\bf Case (1)}: If $\lambda<\omega$ then $F_0(\alpha)=\overline{F_0(H)}$  and by \cref{dimm,lamop,lmono}, we have $\Adim F_0(\alpha)=\Adim F_0(H)=\lambda-1$. The induction hypothesis and \cref{equicat} imply that $F_0(\alpha)\in\ASupp (\cA/\cA_0)_{\lambda-1}=\ASupp F_0(\cA_\lambda)$. Then \cref{bij} implies that $\alpha\in\ASupp\cA_{\lambda}$. 
  {\bf Case (2)}:  If $\lambda \geq \omega$ and $\beta$ is an atom in $\ASpec\cA$ such that $\alpha<\beta$, then $\Adim\beta<\lambda$. Again, by \cref{mondim}, there exists a monoform object $G$ in $\cA$ such that $\overline{G}=\beta$ and $\Adim G=\Adim\beta<\lambda.$ The induction hypothesis implies that $\GKdim G\leq \Adim G$ so that $\GKdim\beta<\lambda$. Hence  $\GKdim\alpha\leq\lambda$ and so $\alpha\in\ASupp\cA_{\lambda}$.  Therefore we have shown that  $\ASupp M\subseteq \ASupp \cA_{\lambda}$ and consequently $\GKdim M\leq\lambda$ by \cref{lemcru}.   
 \end{proof}

  \medskip
  \begin{Corollary}\label{lab}
  Let $\ASpec\cA$ be Alexandrov and $M$ be a locally monoform object in $\cA$. Then $\GKdim\alpha=\Adim\alpha$ for any $\alpha\in\ASupp M$.
  \end{Corollary}
  \begin{proof}
  Given $\alpha\in\ASupp M$, $\Lambda(\alpha)$ is an open subset of $\ASpec\cA$. Then there exits a monoform locally monoform object $H$ such that $\ASupp H=\Lambda(\alpha)$ and $\GKdim H=\GKdim \alpha$. Hence \cref{coindim} implies that $\Adim\alpha=\Adim H=\GKdim \alpha$.
  \end{proof}

   \medskip
   \begin{Corollary}\label{lab1}
  Let $\ASpec\cA$ be Alexandrov and $M$ be a noetherian object in $\cA$. Then there exists $\alpha\in\ASupp M$ such that $\Adim M=\Adim\alpha$. In particular, $\Adim M$ is non-limit ordinal.
  \end{Corollary}
  \begin{proof}
  Since $M$ is noetherian, it is locally monoform and hence $\GKdim M=\Adim M$ by \cref{coindim}. Assume that $\Adim M=\sigma$. We assert that $\sigma$ is a non-limit ordinal. Otherwise, since $M$ is noetherian and $M=\Sigma_{\delta<\sigma}t_{\delta}(M)$, there exists $\rho<\sigma$ such that $M=t_{\rho}(M)$ which is a contradiction. If for any $\alpha\in\ASupp M$, we have $\Adim\alpha<\sigma$, the definition implies that $\Adim M\leq \sigma-1$ which is a contradiction.
  \end{proof}

  \medskip
\begin{Example}\label{rem}
We remark that the equality in the above theorem may not hold if $\ASpec\cA$ is not Alexandrov even if $\cA$ is locally noetherian. To be more precise, if we consider the locally noetherian Grothendieck category$\cA={\rm GrMod} k[x]$ of garded $k[x]$-modules, where $k$ is a field and $x$ is an indeterminate with $\deg x=1$. According to [K1, Example 3.4], $\Adim k[x]=0$ while $\GKdim k[x]=1$. 
\end{Example}


\section{Minimal atoms of an object}
 
In this section we study the minimal atoms of objects in the Grothendieck category $\cA$. Given an object $M$ in $\cA$, an atom $\alpha\in\ASupp M$ is said to be {\it minimal} if it is minimal in $\ASupp M$ under $\leq$. We denote by $\AMin M$, the set of minimal atoms of $M$.
 \medskip

\begin{Lemma}\label{thm}
Let $\sigma$ be an ordinal and let $M$ be an object in $\cA$ with $\Adim M=\sigma$. Then every  $\alpha\in\ASupp M$ with $\Adim \alpha=\sigma$ belongs to $\AMin M$. In particular, if $\ASpec\cA$ is Alexandrov  and $M$ is noetherian, then $\{\alpha\in\ASupp M|\hspace{0.1cm} \Adim\alpha=\sigma\}=\ASpec_{\sigma-1}\cA\cap\ASupp M$ is a finite set.  
\end{Lemma}
\begin{proof}
 If $\alpha\notin\AMin M$, then there exists some $\beta\in\ASupp M$ such that $\beta<\alpha$. But the definition implies that $\Adim\alpha<\Adim \beta$ which contradicts the fact that $\Adim M=\Adim\alpha$. To prove the second claim, it is clear that $M$ is locally monoform and so it follows from \cref{coindim}  that $\GKdim M=\sigma$. Moreover, $\sigma$ is a non-limit ordinal by \cref{lab1}. Since $M$ is noetherian, then $F_{\sigma-1}(M)$ has finite length and so $\ASupp F_{\sigma-1}(M)$ is a finite set.  On the other hand, by \cref{lab} and \cref{p1}, we have 
 $$F_{\sigma-1}(\{\alpha\in\ASupp M|\hspace{0.1cm} \Adim\alpha=\sigma\})=F_{\sigma-1}(\{\alpha\in\ASupp M|\hspace{0.1cm} \GKdim\alpha=\sigma\})$$$$= \ASupp F_{\sigma-1}(M)=F_{\sigma-1}(\ASpec_{\sigma-1}\cA\cap\ASupp M)$$ and hence $\{\alpha\in\ASupp M|\hspace{0.1cm} \Adim\alpha=\sigma\}=\ASpec_{\sigma-1}\cA\cap\ASupp M$ is a finite set by  \cref{bij}.
\end{proof}

The following corollary has been proved by Kanda [K1, Proposition 4.7], where $\cA$ was locally noetherian. Here we prove it without this condition. 
\begin{Corollary}
Let $M$ be a noetherian object in $\cA$ and let $\alpha\in\ASupp M$. Then there exists an atom $\beta$ in $\AMin M$ such that $\beta\leq \alpha$.  
\end{Corollary}
\begin{proof}
Assume that $\alpha\in\ASupp M$ and assume that $F:\cA\To \cA/\cX(\alpha)$ is the canonical functor. We notice that $\ASupp F(M)=F(\ASupp M\cap\overline{\{\alpha\}})$. We observe that $F(M)$ is noetherian.  Assume that $\Adim F(M)=\sigma$. Then by \cref{lab1}, there exists $F(\beta)\in\ASupp F(M)$ such that $\Adim F(\beta)=\sigma$. Hence \cref{thm} implies that $F(\beta)\in\AMin F(M)$. Now, \cref{bij} and \cref{leqf} imply that $\beta\in\AMin M$.
\end{proof}

We now present the first main result of this section which provides a sufficient condition for finiteness of the number of minimal atoms of a noetherian object.
\medskip

\begin{Proposition}\label{finmin}
Let $M$ be a noetherian object in $\cA$. If $\Lambda(\alpha)$ is an open subset of $\ASpec\cA$ for any $\alpha\in\AMin M$, then $\AMin M$ is a finite set. 
\end{Proposition}
\begin{proof}
Let $\alpha\in\AMin M$ and set $W(\alpha)=\{\beta\in\ASpec \cA|\hspace{0.1cm} \alpha<\beta\}$. It is straightforward to show that $W(\alpha)=\Lambda(\alpha)\setminus\overline{\{\alpha\}}$; and hence $W(\alpha)$ is an open subset of $\ASpec\cA$. Consider $\Phi=\cup_{\alpha\in\AMin M}W(\alpha)$, the localizing subcategory $\cX=\ASupp^{-1}(\Phi)$ and the canonical functor $F:\cA\To\cA/\cX$. It follows from [K1, Lemma 5.16] that $\ASupp F(M)=F(\AMin M)$. We notice that for any $\alpha\in\AMin M$, we have $\Lambda(\alpha)\cap(\ASpec\cA\setminus\Phi)=\{\alpha\}$; and hence using \cref{lamop}, $\Lambda(F(\alpha))=\{F(\alpha)\}$ is an open subset of $\ASpec\cA/\cX$ so that $F(\alpha)$ is a maximal atom of $\ASpec\cA/\cX$ by using [S1, Proposition 3.2]. On the other hand, according to [Po, Chap 5, Lemma 8.3], the object $F(M)$ is noetherian. Thus the previous argument implies that $F(M)$ has finite length so that $F(\AMin M)$ is a finite set. Since  $\AMin M\subseteq \ASpec\cA\setminus\ASupp\cX$, the set  $\AMin M$ is finite using \cref{bij}. 
\end{proof}

\medskip

The following example shows that the above result  may not hold in a more general case even if $\cA$ is locally noetherian. 
	
  \begin{Example}([Pa, Example 4.7], [K1, Example 3.4])\label{ex}
It should be noted that the set of minimal atom of a Grothendieck category is not finite when $\cA$ does not have a notherian generator. To be more precise, let $\cA={\rm GrMod} k[x]$ be the category of garded $k[x]$ modules, where $k$ is a field and $x$ is a indeterminate with $\deg x=1$. We notice that $\cA$ is a locally noetherian Grothendieck category. For each $i\in\bbZ$, the object $S_i=x^ik[x]/x^{i+1}k[x]$ is $0$-critical; and hence $\overline{S_i}$ is a minimal atom of $\cA$ for each $i\in\bbZ$. Furthermore, the set of minimal atom of a notherian object is not finite in general even if $\cA$ is locally noetherian. If we consider the noetherian $k[x]$-module $M=k[x]$, then it is easy to see that $\AMin M=\ASupp M=\{\overline{S_j}|\hspace{0.1cm} j\leq 0\}\cup \{\overline{M}\}$. 
\end{Example}

\section{In the case of noncommutative rings}

In this section, we assume that $A$ is a noncommutative ring.  The atom spectrum $\ASpec\textnormal{Mod-}A$ is denoted by $\ASpec A$.  We first recall the classical Krull dimension of right $A$-modules [GW]. 

\begin{Definition}
We define the Krull dimension of right $A$-modules by transfinite induction, introducing classes $\cK_{\sigma}$ of modules, for all ordinals  $\sigma$. Let $\cK_{-1}$ be the class containing precisely of the zero module. Consider an ordinal $\sigma\geq 0$ and suppose that $\cK_{\beta}$ has been defined for all ordinals $\beta<\alpha$. We define $\cK_{\alpha}$, the class of those modules $M$ such that, for every (countable) descending chain $M_0\geq M_1\geq \dots $of submodules of $M$, we have $M_i/M_{i+1}\in\bigcup_{\beta<\alpha}\cK_{\beta}$ for all but finitely many indices $i$. The smallest such $\alpha$ such that $M\in\cK_{\alpha}$ is the {\it Krull dimension} of $M$, denoted by $\Kdim M$ and we say that $\Kdim M$ exists. We notice that the modules of Krull dimension zero are precisely the nonzero artinian modules.   
\end{Definition}

\medskip
\begin{Definition}
 We recall from [Sm1] that a ring $A$ is said to be {\it right fully bounded} if for every prime ideal $\frak p$ of $A$, the ring $A/\frak p$ has the property that every essential right ideal contains a nonzero two-sided ideal. According to [St], a right noetherian right fully bounded ring is called {\it fully right bounded}.
 \end{Definition}
 
 We recall the compressible modules which  have a key role in our studies.

  \begin{Definition} 
We recall from [Sm1] that a non-zero right $A$-module $M$ in $\cA$ is {\it compressible} if each non-zero submodule $L$ of $M$ has some submodule isomorphic to $M$. 
\end{Definition}

\medskip
The following proposition is crucial in our investigation.

\begin{Proposition}[Sm1, Proposition 26.5.10]\label{smp}
Let $A$ be a right fully bounded ring with right Krull dimension. Then every nonzero right $A$-module contains a compressible monoform submodule.
\end{Proposition}

 For a  right fully bounded ring $A$, we investigate when $\ASpec A$ is Alexandrov.
 
\begin{Proposition}\label{two}
Let $A$ be a right fully bounded ring with right Krull dimension. Then $\ASpec A$ is an Alexandrov topological space. 
 \end{Proposition}
\begin{proof}
Assume that $\alpha$ is an atom in $\ASpec A$ and $M$ is a monoform right $A$-module such that $\alpha=\overline{M}$. It follows from \cref{smp} that $M$ contains a compressible monoform submodule $H$ such that $\alpha=\overline{H}$. Therefore $\Lambda(\alpha)=\ASupp H$ by \cref{SaS}. Now the result follows by \cref{alex}.
\end{proof}

\medskip

The following example due to Musson [M] shows that if $A$ is not a right fully bounded ring, then  \cref{two} may not hold. Moreover, \cref{finmin} may not hold even for a cyclic module.

\begin{Example}\label{eexx}
Let $k$ be an algebraically close field of characteristic zero and let $B=k[[t]]$ be the formal power series ring over $k$ in an indeterminate $t$. Define a $k$-linear derivation $\delta$ on $S$ according to the rule $\delta(\Sigma_{n=0}^{\infty}a_nt^n)=\Sigma_{n=0}^{\infty}na_nt^n$. Now, assume that $A=B[\theta]$ is the formal linear differential operator ring (the Ore extension) over $(B,\delta)$. Thus additively, $A$ is the abelian group of all polynomials over $B$ in an indeterminate $\theta$, with a multiplication given by $\theta b=b\theta+\delta(b)$ for all $b\in B$. Since $B$ is noetherian, using [R, Theo$\grave{{\rm r}}$eme 2, p.65], the ring $A$ is right and left noetherian and there is a $B$-isomorphism $B=A/\theta A$. In view of [Go], the nonzero right $A$-submodules of $B$ form a strictly descending chain $B>tB>t^2B>\dots$ and $B$ is a critical right $A$-module of Krull dimension one and so all factors $t^nB/t^{n+1}B$ have Krull dimension zero. Also none of these submodules can embed in any strictly smaller submodule; and hence none of these submodule is compressible. It therefore follows from [GR, Theorem 8.6, Corollary 8.7] that that $A$ is not a right fully bounded ring. Since $k$ is algebraically close field, the maximal two-sided ideals are precisely $\frak m_{\lambda}=(\theta-\lambda)k[\theta]+tA$ with $A/M_{\lambda}\cong k$ for all $\lambda\in k$. Furthermore, for each $n\geq 0$, we have an isomorphism $t^nB/t^{n+1}B\cong A/\frak m_n$ which are pairwise non-isomorphic simple right $A$-modules. Moreover, one can easily show that $\ASupp t^nB=\{\overline{B}\}\cup\{\overline{A/\frak m_i}|\hspace{0.1cm} i\geq n\}$ for every $n\geq 0$; and hence $\{\overline{B}\}=\bigcap_{n\geq 0}\ASupp t^nB$. It now follows from [K1, Proposition 4.4] that $\overline{B}$ is maximal under $\leq$ in $\ASpec A$ so that $\AMin B=\{\overline{B}\}\cup\{\overline{A/\frak m_n}|\hspace{0.1cm} n\geq 0\}$. We also observe that $\ASpec A$ is not Alexandrov as $\{\overline{B}\}$ is not an open subset of $\ASpec\cA$.  
\end{Example}

\medskip
\begin{Definition}
A {\it polynomial identity} on a ring $R$ is a polynomial $p(x_1,\dots,x_n)$ in noncommutative variables $x_1,\dots,x_n$ with coefficients  from $\mathbb{Z}$ such that $p(r_1,\dots,r_n)=0$ for all $r_1,\dots,r_n\in R$. A {\it polynomial identity ring}, or {\it P. I. ring} for short, is a ring which satisfies some monic polynomial identity $p(x_1,\dots,x_n)$ (that is, among the monomials of highest total degree which appear in $p$, at least one has coefficient 1).
\end{Definition}

\begin{Corollary}\label{comend}
Let $R$ be a commutative ring. If  $M$ is a noetherian $R$-module with $A=\End_R(M)$, then $\ASpec A$ is an Alexandrov space.    
\end{Corollary}
\begin{proof}
 Since $M$ is a finitely generated $R$-module, there exists a positive integer $n$ such that $A$ is an $R$-submodule of $M^n$; and hence $A$ is a noetherian $R$-module. Let $I$ be a  right ideal $A$. For every $r\in R$ and $f\in I$ the multiplication $rf=f\circ(r.)$ yields that  $I$ is an $R$-submodule of $A$. Therefore, $A$ is a right noetherian ring; and hence $A$ has right Krull dimension. On the other hand, according to [GW, p. xv], $A$ is a P. I. ring. Hence $A$ is a right fully bounded ring by [Row, Proposition 6.1.48]. Now \cref{two} implies that $\ASpec A$ is Alexandrov.
\end{proof}

\begin{Proposition}\label{irrComp}
Let $\frak p$ be a prime ideal of $A$ and consider the following conditions.\\
\textnormal{(1)} $\frak p$ is an irreducible right ideal.\\
\textnormal{(2)} $A/\frak p$ is a compressible right $A$-module.\\
\textnormal{(3)} $A/\frak p$ is a monoform right $A$-module.

Then the implication \textnormal{(1)}$\Longrightarrow\textnormal{(2)}$ holds if $A$ is a right fully bounded ring, the implication \textnormal{(2)}$\Longrightarrow\textnormal{(3)}$ holds if $A/\frak p$ is locally monoform and \textnormal{(3)}$\Longrightarrow\textnormal{(1)}$ holds.
\end{Proposition}
\begin{proof}
(1)$\Rightarrow$(2). If $\frak p$ is an irreducible right ideal, then every non-zero right submodule of $A/\frak p$ is essential. Given a non-zero submodule $K$ of $A/\frak p$, since $\Ass(K)=\{\frak p\}$, there exists a non-zero element $x\in K$ such that Ann$(xA)=\frak p$. Observe that $\frak p\subseteq \rAnn(x)=\{a\in A|\hspace{0.1cm} xa=0\}$. If $\frak p\subset\rAnn(x)$, since $A$ is fully right bounded, there exists a two-sided ideal $\frak b$ such that $\frak p\subsetneq \frak b\subset \rAnn(x)$. But this implies that $\frak b\subset\Ann(xA)=\frak p$ which is impossible. Thus $\frak p=\rAnn(x)$; and hence $xA\cong A/\frak p$. (2)$\Rightarrow$(3) By the assumption, $A/\frak p$ contains a monoform subobject $K$. Since $A/\frak p$ is compressible, $K$ contains a monoform subobject isomorphic to $A/\frak p$ which implies that $A/\frak p$ is monoform. (3)$\Rightarrow$(1). Since $A/\frak p$ is monoform, any non-zero right submodule is essential. Thus $\frak p$ is an irreducible right ideal.  
\end{proof}

 \medskip
 \begin{Proposition}\label{ircomp}
   Let $A$ be a ring.
  
\textnormal{(1)} If $M$ is a compressible right $A$-module, then $\Ann(M)$ is a prime ideal of $A$.

\textnormal{(2)} Let $A$ be a right fully bounded ring. Then a right $A$-module $M$ is compressible with irreducible right ideal $\Ann(M)$ if and only if $M$ is isomorphic to a submodule of $A/\frak p$ for some prime ideal $\frak p$ that is an irreducible right ideal.
 \end{Proposition} 
\begin{proof}
(1) Setting $\frak p=\Ann(M)$, for any nonzero submodule $K$ of $M$, it is clear that $\frak p\subseteq{\rm Ann}(K)$. On the other hand, since $M$ is compressible, it is isomorphic to a submodule of $K$. Then we have  ${\rm Ann}(K)\subseteq \frak p$; and consequently ${\rm Ann}(K)=\frak p$. For two-sided ideals $\frak a$ and $\frak b$ of $A$ with $\frak a\frak b\subseteq \frak p$, if $K\frak a\neq 0$, the first argument forces that Ann$\frak (K\frak a)=\frak p$ so that $\frak b\subseteq \frak p$. Hence $\frak p$ is a prime ideal of $A$. 

(2) Assume that $M$ is a compressible and set $\frak p=\Ann(M)$. Then  for any nonzero element $x\in M$, a similar argument as mentioned in (1) yields that $\Ann(xA)=\frak p$. Since by the assumption, $\frak p$ is an irreducible right ideal, the proof of (1)$\Longrightarrow$(2) in \cref{irrComp} yields $\frak p=\rAnn(x)$ so that $xA\cong A/\frak p$. Consequently, since $M$ is compressible, it is isomorphic to a submodule of $A/\frak p$. Conversely, it follows from \cref{irrComp} that $A/\frak p$ is a compressible right $A$-module; and hence it is clear that $M$ is a compressible right $A$-module and $\Ann (M)=\frak p$. 
\end{proof}

\medskip
\begin{Corollary}\label{crop}
Let $A$ be a  right fully bounded ring and $\frak p\in\Irr A$ be such that $A/\frak p$ is locally monoform. Then  $\Lambda(\overline{A/\frak p})$ is an open subset of $\ASpec A$.
\end{Corollary}
\begin{proof}
It follows from \cref{irrComp} that $A/\frak p$ is a monoform and compressible right $A$-module. Moreover, \cref{SaS} implies that $\Lambda(\overline{A/\frak p})= \bigcap_{\overline{H}=\overline{A/\frak p}}\ASupp H=\ASupp A/\frak p$ is an open subset of $\ASpec A$.
\end{proof}
\medskip

Let $\alpha\in\ASpec A$ have a noetherian representative $H$. Since $E(\alpha)$ is an indecomposable injective  right $A$-module, $\Ass(H)$ has only one element and we denote it by $\ass(\alpha)$. 

\begin{Corollary}
Let $A$ be a right fully bounded ring and $M$ be a noetherian $A$-module such that $\ass(\alpha)$ is an irreducible right ideal for every $\alpha\in\AMin M$. Then $\AMin M$ is a finite set.  
\end{Corollary}
\begin{proof}
Let $\alpha\in\AMin M$ and $\frak p\in\Ass (\alpha)$. Since $M$ is noetherian, there exists a monoform right noetherian $A$-module $H$ such that $\alpha=\overline{H}$   and $\Ann(H)=\frak p$. We may assume that $H=xA$ for some $x\in H$; and hence $A/\frak p\cong xA$ as $A$ is a right fully bounded. On the other hand, since $A/\frak p$ is a right noetherian $A$-module, it is locally monoform. Hence, \cref{crop} implies that $\Lambda(\alpha)$ is an open subset of $\ASpec A$. Now, the result follows from  \cref{finmin}.
\end{proof}

\medskip

 For a compressible right $A$-module $M$, we denote by $[M]$ the class of all compressible modules $N$ such that $\Ann(M)=\Ann(N)$. Moreover, we set 
\begin{center}
$\Comp A=\{[M]|\hspace{0.05cm}$ $M$ is a compressible right $A$-module$\}$. 
\end{center}

 \medskip
 \begin{Lemma}\label{repe}
 Let $\frak p$ be a prime ideal of $A$ and $M$ be a nonzero compressible right $A$-submodule of $A/\frak p$. Then $\Ann(M)=\frak p$.
 \end{Lemma}
 \begin{proof}
Let  $\Ann(M)=\frak b$.  Without loss of generality we may assume that $M=xA$ for some nonzero element $x$ in $M$. Let $M=\frak a/\frak p$ for some right ideal $\frak a$ of $A$. Then $A\frak a\frak b\subseteq \frak p$ and so $\frak b= \frak p$ as $\frak p$ is a prime ideal of $A$. 
 \end{proof}
Let $A$ be a right fully bounded ring with right Krull dimension. For any prime ideal $\frak p$ of $A$, it follows from \cref{smp} that the right $A$-module $A/\frak p$ contains a compressible monoform submodule $\cm(\frak p)$. If $M$ is another compressible monoform right submodule of $A/\frak p$, then  \cref{repe} implies that $[M]=[\cm(\frak p)]$. We now the following theorem. 
\medskip

\begin{Theorem}\label{tcorr}
Let $A$ be a right fully bounded ring with right Krull dimension. Then there is a bijective map $\Ann:\Comp A\To\Spec A$ given by $[M]\mapsto \Ann(M)$. The inverse map is  $[\cm]:\Spec A\To\Comp A$ given by $\frak p\mapsto [\cm(\frak p)]$.  
\end{Theorem}
\begin{proof}
 Let $M$ be a compressible right $A$-module. Then  $\Ann(M)$ is a prime ideal of $A$ by \cref{ircomp}. Hence, we have  $[\cm(\Ann(M))]=[M]$ by \cref{repe}. On the other hand, it is clear that $\Ann[\cm(\frak p)]=\frak p$ for any prime  ideal $\frak p$ of $A$. 
\end{proof}

\medskip
Let $A$ be a right noetherian ring. For any $\frak p\in\Spec A$, by \cite[Proposition 1.9, Chap VII]{St}, we have $E(A/\frak p)=\bigoplus_i E_i$, where $E_i$ are indecomposable injective modules and isomorphic to each other. Then for each $i$, $\AAss A/\frak p=\AAss E(A/\frak p)=\AAss E_i$ has only one element and we show it by $\aass(\frak p)$.

Let $A$ be a fully right bounded ring. Since $A$ is a right noetherian ring, by [GW, Lemma 15.3], it has right Krull dimension. Then for any atom $\alpha\in\ASpec A$, by \cref{smp}, there exists a compressible monoform right $A$-module $M$ such that $\alpha=\overline{M}$ and $\Ann(M)$ is a prime ideal of $A$ by \cref{ircomp}. If $N$ is another such compressible monoforom right $A$-module, then $\Ann(M)=\Ann(N)$ as $M$ and $N$ have  a common nonzero submodule. In this case, we denote $\Ann(M)$ by $\Ann(\alpha)$. We now the following theorem.

\medskip
 \begin{Theorem}\label{bijj}
Let $A$ be a  fully right bounded ring. Then there is a bijective map $\aass:\Spec A\To\ASpec A$ given by $\frak p\mapsto \aass(\frak p)$.  The inverse map is $\Ann:\ASpec A\To\Spec A$ given by $\alpha\mapsto \Ann(\alpha)$. 
\end{Theorem}
\begin{proof}
 For any $\alpha\in\ASpec A$, there exists a compressible monoform right $A$-module $H$ such that $\alpha=\overline{H}$. Then $\Ann(\alpha)=\Ann(H
)=\frak q_1$ is a prime ideal of $A$. We show that $\aass(\Ann(\alpha))=\alpha$. If $\aass(\Ann(\alpha))=\beta$, then there exists a compressible monoform right $A$-module $G$ such that $E(G)=E(\beta)$ is a submodule of $E(A/\frak q_1)$. If $\ass E(\alpha)=\frak q_2$, then $E(\alpha)=E(H)$ is a submodule of $E(A/\frak q_2)$; and hence $\frak q_1=\frak q_2$ by \cref{repe}. Since $A$ is fully right bounded, we have $E(\alpha)\cong E(\beta)$ so that $\alpha=\beta$. On the other hand, for any prime ideal $\frak p$ of $A$, it follows from \cref{repe} that $\Ann(\aass(\frak p))=\frak p$.
\end{proof}

The above theorems yields the following corollary.  
\begin{Corollary}\label{cabij}
Let $A$ be a fully right bounded ring. Then there is a bijective map $\Comp A\To \ASpec A$, given by $[M]\mapsto\overline{\cm(\Ann(M))}$. The inverse map is given by $\alpha\mapsto [\cm(\Ann(\alpha)]$. 
\end{Corollary}
\medskip
\begin{Definition and Remark}Let  $A$ be a fully right bounded ring. For any right $A$-module $M$, by the above corollary, we define the {\it compressible support} of $M$ that is 
$$\CSupp M=\{[X]\in\Comp A|\hspace{0.1cm} \overline{\cm(\Ann(X))}\in\ASupp M\}.$$
For any exact sequence $0\To N\To M\To K\To 0$ of right $A$-modules, it is clear that $$\CSupp M=\CSupp K\cup \CSupp L.$$
By \cref{cabij} and \cite[Proposition 3.2]{K2}, the collection $\{\CSupp M|\hspace{0.1cm} M$ is a right $A$-module$\}$ forms a basis  of open subsets of a topology on $\Comp A$. We define a preorder relation on $\Comp A$ as follows: For any $[X],[Y]\in\Comp A$ we define
\begin{center}
$[X]<[Y]$ if and only if $\overline{\cm(\Ann(X))}<\overline{\cm(\Ann(Y))}.$
\end{center}

\end{Definition and Remark}

\medskip
\begin{Definition}Let $A$ be a right noetherian. For any right $A$-module $M$, {\it the  support} $M$  is $$\Supp M=\{\frak p\in\Spec A|\hspace{0.1cm} \aass\frak p\in\ASupp M\}.$$
For any exact sequence $0\To N\To M\To K\To 0$ of right $A$-modules, it is clear that $$\Supp M=\Supp K\cup \Supp L.$$ Also, it should be noted that this notion coincides with the classical support of modules over a commutative ring. 
For a well-ordered subset $W\subseteq \Spec A$, the length $\lambda(W)$ is defined by
 
$$\lambda(W)=\begin{cases}
{\rm ord}\hspace{0.1cm} W \hspace{0.1cm} {\rm if \hspace{0.1cm}the\hspace{0.1cm} order\hspace{0.1cm} type \hspace{0.1cm} of}\hspace{0.1cm} W\hspace{0.1cm} {\rm is\hspace{0.1cm} a\hspace{0.1cm} limit\hspace{0.1cm} ordinal};\\
  {\rm ord}\hspace{0.1cm}W-1 \hspace{0.1cm} {\rm otherwise}.
\end{cases}$$

By [GR], the {\it little  Krull dimension} of a right $A$-module $M$ is denoted by $\kdim M$ that is 
 \begin{center}
  $\kdim M=\sup\{\lambda(C)|\hspace{0.1cm} C$ is a chain of prime ideals in $\Supp M\}$. 
 \end{center}
  
Similarly, for a fully right bounded ring $A$, the {\it compressible Krull dimension} of $M$ is dented by $\ckdim M$ that is 
 \begin{center}
 $\ckdim M=\sup\{\lambda(C)|\hspace{0.1cm} C$ is a chain of elements in $\CSupp M\}$.
\end{center}
\end{Definition}
\medskip

\begin{Lemma}\label{ineqq}
Let $A$ be a fully right bounded ring and $\frak p\subset\frak q$  be primes ideal of $A$. Then we have  $\aass\frak p<\aass\frak q.$
\end{Lemma}
\begin{proof}
 By \cref{smp}, there exist compressible monoform right $A$-modules $H$ and $G$ such that $\Lambda(\aass(\frak p))=\ASupp H=\ASupp A/\frak p$ and $\Lambda(\aass(\frak q))=\ASupp G=\ASupp A/\frak q$. Therefore the inclusion $\frak p\subset \frak q$ implies that $\Lambda(\aass(\frak q))\subset \Lambda(\aass(\frak p))$; and consequently $\aass\frak p<\aass\frak q$. 
\end{proof}

\medskip
 \begin{Corollary}\label{bkr}
Let $A$ be a  fully right bounded ring and $M$ be a right $A$-module with $\Adim M<\omega$. Then $\GKdim M=\Adim M=\kdim M=\ckdim M$.
\end{Corollary}
\begin{proof}
 Let $\Adim M=n$ and $\cA=\Mod\cA$. Since $A$ is right noetherian, $M$ is locally monoform and $A$ has right Krull dimension by \cite[Lemma 15.3]{GW}. Then $\ASpec A$ is Alexandrov by \cref{two}. To prove the first equality, it suffices to prove that $\GKdim M\leq \Adim M$ and so the result follows by \cref{coindim}.  We prove  by induction on $n$. If $n=0$, then  any atom  $\alpha\in\ASupp M$  is maximal under $\leq$  and so $\alpha$ is maximal by \cref{alexmax}. Therefore $\alpha\in\ASupp\cA_0$.  This implies that $\ASupp M\subset\ASupp\cA_0$ and so $\GKdim M=0$ by \cref{lemcru}. Now assume that $n>0$. By \cref{mondim}, there exists a monoform object $H$ in $\cA$ such that $\alpha=\overline{H}$ and $\Adim\alpha=\Adim H$. If $\Adim\alpha<n$, the induction hypothesis implies that $\GKdim H<n$; and hence  $\alpha\in\ASupp\cA_{n}$. If $\dim \alpha=n$, then $F_0(\alpha)=\overline{F_0(H)}$  and by \cref{dimm,lamop,lmono}, we have $\Adim F_0(\alpha)=\Adim F_0(H)=n-1$. The induction hypothesis and \cref{equicat} imply that $F_0(\alpha)\in\ASupp (\cA/\cA_0)_{n-1}=\ASupp F_0(\cA_n)$. Hence \cref{bij} implies that $\alpha\in\ASupp\cA_{n}$; and consequently $\ASupp M\subset \ASupp\cA_{n}$. Now, \cref{lemcru} implies that $\GKdim M\leq n$. To prove the second equality, since $\Adim M=n$, there exists a chain  $\alpha_0<\alpha_1<\dots <\alpha_n$ of atoms in $\ASupp M$. Then according to \cref{bijj}, it is clear that $\Ann(\alpha_0)\subset\Ann(\alpha_1)\subset\dots \subset\Ann(\alpha_n)$ is a chain of prime ideals in $\Supp M$ so that $\kdim M\geq n$. If $\kdim M>n$, there exists a chain of prime ideals $\frak p_0\subset \frak p_1\subset \dots \subset \frak p_{n+1}$ in $\Supp M$. Hence, \cref{ineqq} yields a chain of atoms $\aass(\frak p_0)<\aass(\frak p_1)<\dots<\aass(\frak p_{n+1})$ in $\ASupp M$. But this implies that $\Adim M\geq n+1$ which is a contradiction. The equality $\Adim M=\ckdim M$ follows from \cref{cabij}.
\end{proof}

The following proposition shows that the Krull-Gabriel dimension of right $A$-modules serves as a lower bound for the classical Krull dimension as defined above. 
\medskip

\begin{Proposition}\label{kkl}
Let $M$ be a right $A$-module with $\Kdim M=\sigma$. Then $$\GKdim M\leq \begin{cases}
\sigma\hspace{0.5cm}{\rm if} \hspace{0.1cm}\sigma< \omega\\
  \sigma+1 
 \hspace{0.5cm}{\rm if} \hspace{0.1cm}\sigma\geq \omega.
 \end{cases}$$

 In particular, if $M$ is noetherian, the inequalities are equalities.
\end{Proposition}
\begin{proof}
We proceed by induction on $\sigma$. We first consider $\sigma<\omega$. If $\sigma=0$, then $M$ is artinian  and so $\GKdim M=0$. If $\sigma>1$ and $\GKdim M\nleq \sigma$, we have $M\notin\cA_{\sigma}$ and so $F_{\sigma-1}(M)$ is not artinian. Then there exists an unstable descending chain $M'_0\supseteq M'_1\dots$ of submodules of $F_{\sigma-1}(M)$. According to [Po, Chap 4, Corollary 3.10], there exists a descending chain $M_0\supseteq M_1\dots$ of submodules of $M$ such that $F_{\sigma-1}(M_i)=M'_i$ for each $i$ and since $F_{\sigma-1}(M_i/M_{i+1})\neq 0$ for infinitely many indexes $i$, the induction hypothesis implies that  $M_i/M_{i+1}\notin\cK_{\sigma-1}$ for infinitely many indices $i$ which is a contradiction. To prove the second assertion, assume that $M$ is noetherian and so by \cref{coindim} and \cref{lab1}, there exists a non-limit ordinal $\delta$ such that $\GKdim M=\delta$. We proceed by induction on $\delta$ that $\Kdim M\leq \GKdim M$. If $\delta=0$, then $M$ has finite length and so $\Kdim M=0$. If $\delta>1$, since $F_{\delta-1}(M)$ has finite length, for any descending chain  $M_0\supseteq M_1\dots$ of submodules of $M$, there exists some non-negative integer $n$ such that $F_{\delta-1}(M_i/M_{i-1})=0$ for all $i\geq n$ and so the induction hypothesis implies that $\Kdim (M_i/M_{i-1})\leq \delta-1$ so that $\Kdim M\leq \delta$. We now assume that $\sigma\geq \omega$. Then for any descending chain  $M_0\supseteq M_1\dots$ of submodules of $M$, there exists some non-negative integer $n$ such that $\Kdim(M_i/M_{i-1})<\sigma$ for all $i\geq n$. Hence $F_{\sigma}(M_i/M_{i-1})=0$ for all $i\geq n$  by induction hypothesis. This implies that $F_{\sigma}(M)$ is artinian and so $\GKdim M\leq {\sigma+1}$ as $M\in\cA_{\sigma+1}$. If $M$ is noetherian and $\GKdim M=\delta$, we prove by transfinite induction on $\delta$ that $\Kdim M+1\leq \delta$. If $\delta=\omega+1$, the $F_{\omega}(M)$ has finite length and so for any descending chain  $M_0\supseteq M_1\dots$ of submodules of $M$, there exists some non-negative integer $n$ such that $F_{\omega}(M_i/M_{i-1})=0$ for all $i\geq n$ so that $\GKdim (M_i/ M_{i-1})\leq \omega$ for all $i\geq n$. Since the Krull-Gabriel dimension of noetherian modules are non-limit ordinals, using the first case we deduce that $\Kdim (M_i/ M_{i-1})=\GKdim (M_i/ M_{i-1})<\omega$ for all $i\geq n$; and hence $\Kdim M\leq \omega$. If $\delta>\omega+1$, similar to the induction step, $F_{\delta-1}(M)$ has finite length and so  for any descending chain  $M_0\supseteq M_1\dots$ of submodules of $M$, there exists some non-negative integer $n$ such that $F_{\delta-1}(M_i/M_{i-1})=0$ for all $i\geq n$ so that $\GKdim (M_i/ M_{i-1})\leq \delta-1$ for all $i\geq n$. Now, the induction hypothesis implies that $\Kdim (M_i/ M_{i-1})=\GKdim (M_i/ M_{i-1})-1<\delta -1$ for all $i\geq n$; and hence $\Kdim M\leq \delta-1$. 
\end{proof}

\begin{Corollary}
Let $A$ be a right noetherian ring and $\sigma$ be an ordinal such that $\cK_{\sigma}$ is a localizing subcategory of $\Mod A$. Then  $$\cK_{\sigma}= \begin{cases}
\cA_\sigma\hspace{0.5cm}{\rm if} \hspace{0.1cm}\sigma< \omega\\
  \cA_{\sigma+1}
 \hspace{0.5cm}{\rm if} \hspace{0.1cm}\sigma\geq \omega.
 \end{cases}$$
\end{Corollary}
\begin{proof}
The inclusion $\subseteq $ is clear by \cref{kkl}. For the converse, assume that $\sigma<\omega$ and the case $\sigma>\omega+1$ is similar. Let  $M\in\cA_\sigma$. Then $M=\underset{\rightarrow}{\lim}M_i$ is the direct limit of its noetherian submodules. By \cref{kkl}, each $M_i\in\cK_\sigma$ and hence $M\in\cK_\sigma$ as $K_\sigma$ is localizing.
\end{proof}
\medskip

\begin{Corollary}\label{last}
Let $A$ be a fully  right bounded ring and $M$ be a right $A$-module with $\Kdim M<\omega$. Then $\kdim M\leq \Kdim M$. In particular, if $M$ is noetherian, then $\kdim M=\Kdim M$
\end{Corollary}
\begin{proof}
The result follows from \cref{bkr} and \cref{kkl}. 
\end{proof}
 
\medskip
 \begin{Remark}
 We notice that the modules of Krull-Gabriel  dimension zero are precisely the nonzero semi-artinian modules; see for example \cite[Chap VIII]{St}. As semi-artinian modules are not necessarily  artinian, the  inequality in \cref{last} may be strict. To be more precise, if $P$ is the set of prime integers and $M=\bigoplus_{p\in P}\mathbb{Z}/p\mathbb{Z}$, then $\kdim M=0$ by \cref{bkr} while $\Kdim M>0$ as $M$ is not artinian $\mathbb{Z}$-module. 
 \end{Remark}

By [St, Chap VII], a right noetherian ring is said to be {\it right classical} if for any indecomposable injective  right $A$-module $E$  with the associated prime ideal $\frak p$, we have $E=\Gamma_{\frak p}(E)=\bigcup_{n=1}^{\infty}(0:_E\frak p^n)$. It is known that every commutative noetherian ring is classical.

\medskip

\begin{Proposition}\label{laspr}
Let $A$ be a right classical ring and $E$ be an indecomposable  injective right $A$-module with $\ass E=\frak p$. If $E$ has Krull dimension, then $\Kdim E=\Kdim A/\frak p$.
\end{Proposition}
\begin{proof}
By \cite[Chap V, Proposition 5.9, Corollary 5.10]{St} and \cite[Chap VII, Proposition 1.9]{St}, we have $E(A/\frak p)=\bigoplus_{i=1}^nE_i$ such that all $E_i$ are indecomposable and isomorphic to each other. On the other hand since $E$ is a direct summand of $E(A/\frak p)$, \cite[Chap V, Corollary 5.5]{St} implies that $E$ is isomorphic to $E_i$ for each $i$; and hence $\Kdim E=\Kdim E(A/\frak p)$. Therefore $\Kdim A/\frak p\leq \Kdim E$. We observe that  $(0:_E\frak p^n)$ has Krull dimension and so we show that $\Kdim (0:_E\frak p^n)\leq \Kdim A/\frak p$ for all $n>0$. For $n=1$, since $(0:_E\frak p)=\bigcup L_i$ is the direct union of its finitely generated  right $A/\frak p$-submodules $L_i$  and $\Kdim L_i\leq \Kdim A/\frak p$, it follows from \cite[Chap 6, Lemma 2.17]{MR} that $\Kdim (0:_E\frak p)\leq \dim A/\frak p$. For $n>0$, there exists the following exact sequence  $$0\To (0:_E\frak p^n)\To (0:_E\frak p^{n+1})\To (0:_E\frak p^{n+1})/(0:_E\frak p^n)\To 0.$$  
We notice that  $(0:_E\frak p^{n+1})/(0:_E\frak p^n)$ is an $A/\frak p$-module and it has Krull dimension by the above exact sequence. Then a similar argument as for the case $n=1$ yields that $\Kdim (0:_E\frak p^{n+1})/(0:_E\frak p^n)\leq \Kdim A/\frak p$; and hence $\Kdim (0:_E\frak p^{n+1})\leq \Kdim A/\frak p$. It now follows from \cite[Chap 6, Lemma 2.17]{MR} that $\Kdim E\leq \Kdim A/\frak p$. 
\end{proof}
 The following corollary shows that the equality in \cref{last} may hold without the noetherian  condition. 
\medskip
\begin{Corollary}
Let $A$ be a commutative noetherian ring and  $E$ be an injective $A$-module with Krull dimension. Then $\Kdim E=\kdim E\leq 1.$  
\end{Corollary}
\begin{proof}
Since $E$ is the direct sum of indecomposable injective modules, by \cite[Chap 6, Lemma 2.17]{MR}, there exists a prime ideal $\frak p$ of $A$ such that $\Kdim E=\Kdim E(A/\frak p)$. It follows from \cite[Theorem 2.10]{Sm2}, \cref{laspr} and \cref{last} that 
$$\Kdim E(A/\frak p)=\Kdim A/\frak p=\kdim A/\frak p\leq \kdim E(A/\frak p)\leq 1.$$ 
Hence  \cref{last} implies that $\Kdim E(A/\frak p)=\kdim E(A/\frak p)\leq\kdim E$ and so the result follows.
\end{proof}



\begin{thebibliography}{Saz}
\bibitem[D]{D} J. Dixmier, {\em Alg$\grave{e}$bres Enveloppantes}, Gauthier-Villars-Bruxells-Moteral, 1974.

\bibitem[G]{G} P. Gabriel, {\em Des cat\'egories ab\'elinnes}. Bull. Soc. Math. France {\bf 90} (1962), 323-448.


\bibitem[GR]{GR}
Robert Gordon and J.~C. Robson, \emph{Krull dimension}, American Mathematical Society, Providence, R.I., 1973, Memoirs of the American Mathematical Society, No.~133. 

\bibitem[Go]{Go}
K.~R. Goodearl, \emph{Incompressible critical modules}, Comm. Algebra, {\bf 8} (1980), no.19, 1845-1851.

\bibitem[GW]{GW}
K.~R. Goodearl and R.~B. Warfield, Jr., \emph{An introduction to noncommutative {N}oetherian rings}, second ed., London Mathematical Society Student Texts, vol.~61, Cambridge University Press, Cambridge, 2004. 

\bibitem[K1]{K2} R. Kanda, {\em Specialization orders on atom spectra
of Grothendieck categories}, J. Pure and Appl. Algebra {\bf 219} (2015), 4907-4952.
 
\bibitem[K2]{K3} R. Kanda, {\em  Finiteness of the number of minimal atoms in Grothendieck categories}, J. Algrbra, {\bf 527} (2019), 182-195.


\bibitem[M]{M} I. M. Musson, {\em Some examples of modules over noetherian rings}, Glasgow Math. J, {\bf 23} (1982), 9-13.

\bibitem[MR]{MR} J. C. McConnel and J. C. Robson, {\em Noncommutative noetherian rings}, Wiley-Interscience, New York, 1987.

\bibitem[Pa]{Pa}
C. J. Pappacena, \emph{The injective spectrum of a noncommutative space}, J. Algebra \textbf{250} (2002), no.~2, 559--602. 

\bibitem[Po]{Po} N. Popoescu, {\em Abelian categories with
applications to rings and modules}, London Mathematical Society
Monographs, No. 3, Academic Press, London-New York, 1973.

\bibitem[R]{R} G. Renault, \emph{Alg$\grave{{\rm e}}$bre non Commutative}, Paris(1975), Gauthiers-Villars.

\bibitem[Ro]{Ro}  A. L. Rosenberg, {\em Noncommutative algebraic geometry and representation of quantized algebras}, Kluwer Academic Publishers, 1995.

\bibitem[Row]{Row} Louis H. Rowen, {\em Ring theory: Volume II}, Academic Press, 1988.

\bibitem[S1]{Sa1} R. Sazeedeh, {\em Monoform objects and localization theory in abelian
categories}, J. Homotopy Relat. Struct, {\bf 13}(2018), 443-460.

\bibitem[S2]{S2} R. Sazeedeh, {Classifying localizing subcategories of a Grothendieck category}, arXive: 2507.13749v2 [math. CT] 25 Nov 2025.

\bibitem[SS]{SaS} F. Savoji and R. Sazeedeh, {\em Local Cohomology in Grothendieck categoties}, J. Algebra. Appl, {\bf 19} (2020), no.11, 2050222.

\bibitem[Sm1] {Sm1} P. F. Smith, \emph{Compressible and related modules}, Abelian Groups, Rings, Modules, and homological algebra, 2006, Lecture Notes in Pure and Applied Mathematics, 295-312.

\bibitem[Sm2] {Sm2} P. F. Smith, \emph{Krull dimension of injective modules over commutative noetherian rinfs}, Canad. Math. Bull, {\bf 48} (2005), no. 2, 275-282.

\bibitem[St]{St} B. Stenstrom, {\em Rings of Quotients: An introduction to methods of ring theory}, Springer-Verlag 1975.
 
\end{thebibliography}
\end{document}